# Insight into China's Economically Motivated Adulteration Risk in Online Raw Agricultural Product Sales

Hengyu Liu
Department of Management Science and Engineering
School of Economics and Management
Beijing University of Posts and Telecommunications, Beijing, 100876, China

Wen Tong[*]
Department of Logistics Economics and Statistics
School of Economics and Management
Shanghai Maritime University, Shanghai, 200135, China

**Abstract**

Uncertainty in quality and the inspectors' imperfect testing capability leave raw agricultural products (e.g., fresh produce, seafood, livestock and poultry products, etc.) wide open to economically motivated adulteration (EMA), and the strong demand for online shopping of these products in China makes this situation even worse. In this paper, we develop a game-theoretic framework to investigate online raw agricultural product sellers' *preemptive* EMA behavior on an ecommerce platform (EP). Particularly, the sellers differ from each other in the original quality of their products. We characterize the sellers' equilibrium pricing and adulteration decisions and the EP's optimal take rate decision, and analyze how the sampling inspections and adulteration penalty jointly impact these decisions. Moreover, we investigate three managerial levers, such as claiming a higher-than-law-requires penalty, that the administrative departments or the EP can use to deter EMA. Finally, we use the real-word data from *Taobao.com* to calibrate our model and derive more managerial insights from the analytical findings. We find that the heterogenous sellers' adulteration decisions are symmetric and their ex-post pricing decisions lead them to evenly share the market on the EP. Interestingly, we show that the EP's higher take rate will inhibit the sellers' adulteration behavior. However, the profit-maximizing EP may indulge the sellers' adulteration behavior by intentionally decreasing this rate. Our results highlight a penalty-inspection-centered approach as essential to combat EMA, and the three levels can play a role as supplements to this approach under certain conditions.

---

[*] Corresponding author, email: tongwen@shmtu.edu.cn

# 1. Introduction

Before COVID-19 forever changed retail, many people in China were already buying food online thanks to the country's sophisticated omnichannel retail ecosystem across urbanities. The pandemic provided further impetus for people not to dine out, but to cook their meals at home. However, the significant growth in demand for online food also gives rise to a potential avenue for food safety risks. Statistics show that courts across China heard approximately 49,000 cases related to online shopping contracts from 2017 to the end of June 2020, and more than 45 percent of them involved food; in particular, food safety is a major complaint[1]. Among various food safety risks, a growing concern is the introduction of hazards by deliberate human actions known as ***economically motivated adulteration*** (EMA) (or *food fraud*, Zhang and Xue, 2016). EMA is defined as "*intentional substitution or addition of a substance in a product for the purpose of increasing the apparent value of product or reducing the cost of its production*" [2]. EMA is not just an economic issue, it poses a serious threat to public health. For example, the notorious 2008 "melamine adulteration scandal" revealed that farmers and collectors intentionally added melamine to infant formula to increase the protein content of the milk, causing over 300,000 illnesses, 50,000 hospitalizations, and at least 6 deaths. Although the majority of the publicly known EMA incidents originated from developing countries, they have also been a major concern for developed countries (Levi et al., 2020), e.g., the 2008 "pine nuts incident" in the US and the 2013 "horse meat sandal" in the EU. Experts estimate that EMA affects 1% of the global food industry at a cost of about $10 -$15 billion a year[3].

The situation of EMA in raw agricultural products is even worse than in processed foods. This is because, unlike processed food that is produced under uniform quality standards, there usually exists quality uncertainty and thus quality heterogeneity in raw agricultural products, even if they are obtained from the same farmland, leaving them wide open to an EMA hazard. In practice, quality uncertainty can be a major cause of EMA in markets with quality-based pricing. Specifically, EMA can occur *before* the uncertainty is resolved, i.e., ***preemptive* EMA**, which seeks to decrease the likelihood of deriving low-quality products. *Preemptive* EMA has been a serious concern in pork, poultry and seafood farming in various countries, including China, India, Bangladesh, and Vietnam (Doyle et al., 2013; Levi et al., 2020). In these countries, some unethical farms may overuse antibiotics to prevent producing sick or underweight animals. For example, the 2012 KFC "instant chicken scandal" in China revealed that the chickens used by KFC were treated with 18 illegal antibiotics to avoid avian flu[4]. More recently, China's 2020 "CCTV 3.15 Gala" reported that some sea cucumber farmers put dichlorvos (an insecticide) into the pond to remove the organisms that are not conducive to sea cucumber growth[5].

Another factor that contributes to EMA in raw agricultural products is the ***imperfect***

***quality testing capability*** of inspectors, including the administrative departments (e.g., China's State Administration for Market Regulation, CSAMR) and the professional testing labs. The implications of "imperfect" are twofold. The first is the *limited testing coverage*. In China, inspectors implement food quality surveillance and monitoring mainly via market access control and spot checks[6]. However, like many other developing countries where highly dispersed smallholder farmers make up the majority of the farmer population, it is nearly impossible to inspect every farmer or trace every unit of supply back to the producing source (Levi et al., 2020). The second is the *limited testing accuracy*. Currently, unethical farmers can engineer fraudulent ingredients to evade existing quality assurance (QA) and quality control (QC) systems. In fact, adulterants are increasingly unconventional, and current food protection systems are unable to identify the nearly infinite number of potential adulterants that may show up in the food supply (Moore et al., 2011, 2012). Taking melamine as an example, before the 2008 "melamine adulteration scandal", this ingredient was considered neither a potential contaminant, nor an adulterant in routine QA/QC analyses (Spink, 2011; Moore et al., 2011).

Adding to the inspectors' imperfect testing capability, online selling fuels EMA in raw agricultural products. On the one hand, e-commerce platforms (EPs) (e.g., *Taobao.com*, *JD.com*, and *Pinduoduo*) substantially increase farmers' access to broader markets, which was especially true during the pandemic. In 2020, *Alibaba* launched the *Rural Support Program* to help China's farmers connect with its 41 million followers, and 15 million kilograms of fresh produce were sold during the first three days[7]. It was estimated that approximately 610 billion yuan-worth (1$≈6.32 yuan) agricultural products were sold via EPs in China in 2020, and the scale was expected to reach 789 billion by the end of 2021[8]. Unfortunately, the farmers' easier access to markets also indicates greater EMA risk, as EPs make it less possible for inspectors to achieve full coverage and on-site inspections of the numerous products sold on them. Taking *JD.com* for example, the fresh produce SKU reaches 10,000[9]. On the other hand, online selling increases information asymmetry between sellers and consumers. This is because consumers have to deal with many faceless sellers who provide different prices and nonuniform quality descriptions of their products, while no physical samples are available to them. Consequently, it is difficult for consumers to identify their ideal products. Moreover, even if they receive the adulterated products, they usually cannot detect sellers' misconduct due to the lack of professional knowledge and testing equipment.

Although both the quality characteristics of raw agricultural products and the inspectors' imperfect testing capability facilitate EMA on EPs, when making adulteration decisions, a seller (he) must consider the following three factors. First, the ***adulterated quality*** of his products compared to those of the others on the same EP. When buying raw agricultural products online, a consumer (she) will comprehensively compare the quality and prices claimed by different sellers and choose the ones that bring her the highest utility. Therefore, if a seller can obtain (fake) high-quality products, he can ensure a competitive advantage of his products on the EP. The second is the ***detection risk*** of his adulteration behavior by inspectors (and even consumers). Normally, the detection risk depends on the relative amount of adulterants added to the products because a large amount may change the characteristics of products and even cause adverse symptoms after consumption (Levi et al., 2020). Once caught adulterating, in addition to the direct monetary fines by administrative departments, EPs will remove the unethical sellers' products and even suspend their licenses[10,11]. The third is the ***negative externality*** of his and other sellers' adulteration behavior. China's recent *Measures for the Investigation and Punishment of Illegal Acts Related to Online Food Safety* (2021 Amendment) requires the inspectors to publish the inspection results to the public [12]. For example, CSAMR will periodically publish the spot check results of agricultural products sold on different EPs, including the sellers' and EPs' names and the detailed product information[13]. In this case, if any seller is reported selling adulterated products on the EP, the quality-conscious consumers will reduce consumption on it and are very likely to share the unpleasant information with their friends or family. As a result, the demand on the EP will globally decrease due to the word-of-mouth effect. Taken together, the above three factors will impact sellers' adulteration behavior in opposite ways. Specifically, the first factor will motivate them to adulterate more, while the last two factors will depress their such behavior.

To date, a thorough understanding of raw agricultural product sellers' strategic adulteration behavior on EPs and its implications for different stakeholders is absent in the literature, especially in operations management (OM). Witnessing the broad market prospect of online raw agricultural product sales, it is important to understand and eliminate EMA risk on EPs, which is relevant to public health. This research answers two main research questions: (1) What are the raw agricultural product sellers' equilibrium adulteration and pricing decisions on EP? (2) How do the detection risk, penalty, and the EP's take rate impact their adulteration behavior? In addition, we investigate three managerial levers, such as adopting a traceability system, that the administrative departments or the EP can use to mitigate the EMA risk. The managerial insights derived in this work can not only offer practical guidelines to policy-makers and different stakeholders on EPs in China, but can also provide recommendations for stakeholders in other countries that are suffering from the EMA risk in online selling.

To address these questions, we consider a network of sellers (including farmers) sell certain kind of raw agricultural products on an EP, and they differ from each other in the original quality of their products. To increase the selling prices and orders, the sellers may engage in EMA considering the detection risk and the negative externality of their adulteration behavior, both of which are increasing in the amount of adulterants added. Specifically, we investigate the *preemptive* EMA scenario under imperfect testing in this paper. Moreover, the EP will decide the take rate charged from its sellers to maximize its expected profit. By modeling the two stakeholders' decisions a three-stage Stackelberg game, we first analyze the sellers' equilibrium selling prices for their adulterated products by using the *quality-price-ratio*-based (QPR-based) choice model under the QPR maximization consumer choice criterion (Xie et al., 2021). Second, we derive their equilibrium adulteration decisions under the *preemptive* EMA scenario. Particularly, when making these decisions, quality uncertainty is yet to be realized. Third, we derive the EP's optimal take rate decision taking into account its impact on the sellers' future decisions. After that, we derive selles' equilibrium adulteration decisions under the *reactive* EMA scenario. When making these decisions, quality uncertainty has been realized. Then we analyze three managerial levers that enable the administrative departments (i.e., setting an administrative penalty to avoid EP's malpractice regarding food safety management) or the EP (i.e., claiming a higher-than-law requires adulteration penalty to its consumers and adopting a traceability system) to deter EMA risk. Finally, we calibrate our model with data from *Taobao.com*, China's largest online shopping platform, and more managerial insights from the analytical findings. Our analysis generates the following managerial insights:

(1) Interestingly, the heterogeneous sellers will make symmetric adulteration decisions on the EP, since they share the negative externality of each other's misconduct and their pricing decisions are expost. Moreover, this decision is decreasing in the *overall penalty risk*, which is positively affected by the proportion of sampling inspections and adulteration penalties. More specifically, if this risk is sufficiently high (low), they will give up adulterating (adulterate to their maximum levels); otherwise, they will adulterate to some extent, but not to the maximum level. Furthermore, their pricing decisions are to make the QPRs of their adulterated products as competitive as the others, and they will evenly share the market on the EP. These findings indicate that a penalty-inspection-centered approach that ensures a sufficiently high *overall penalty risk* on the EP is essential to eliminate the sellers' adulteration behavior.

(2) Contrary to our intuition that the sellers will adulterate more if the EP increases the take rate, this higher rate will inhibit their adulteration behavior. To this end, how the take rate impacts the EP's expected profit is nonmonotonic and is also determined by the *overall penalty risk* faced by its sellers. We characterize the EP's optimal take rate decision and demonstrate that such a decision will very likely induce the sellers' severe adulteration behavior. Specifically, the *zero-adulterant* scene can only be achieved under the case of a very high *overall penalty risk* on the EP. This unintended implication suggests that a better-designed penalty-inspection-centered approach needs to be implemented together with the requirement of the EP's large enough take rate.

(3) If the proportion of sampling inspections from the administrative departments' inspectors is relatively small or moderate, an unreliable EP has the motive to cover up the sellers' adulteration behavior, and such motive is even stronger as its own inspectors inspect more sellers. In this case, setting a large administrative penalty can effectively avoid the EP's such malpractice. Moreover, if the EP's take rate is limited and the adulteration penalty is moderate, it has the incentive to claim a higher-than-law-requires penalty to its consumers. Furthermore, when the R&D expense of the traceability system is small enough and the sellers face a moderate *overall penalty risk* on the EP, the two stakeholders can achieve a win–win situation by adopting this system.

We organize the rest of our paper as follows: In Section 2, we review the relevant literature and identify our research gaps. In Section 3, we describe the model formulations. In Section 4, we first introduce the QPR-based consumer choice model, and then derive the sellers' equilibrium pricing and adulteration decisions as well as the EP's optimal take rate decision. In Section 5, we exploit the value of three managerial levers that can be used to reduce the EMA risk. In Section 6, we conduct a numerical study using data from *Taobao.com* and explore more managerial insights. All the proofs are postponed to Appendix B.

## 2. Literature review

The first stream of literature related to this study is on food safety and quality management. We refer readers to Gorris (2005), Nagurney and Li (2016), and Chen et al. (2021) for comprehensive reviews in this area. Here, we only review the works published in the last three years. Fan et al. (2020) developed a Stackelberg game model to derive the subgame perfect equilibria strategies for a 2-tier supply chain. They found that the liability cost of low-quality product sharing could affect pricing decisions and product quality. Zhao et al. (2021) found that internal and external supply chain integration are the key factors to improve product quality. Ayvaz-cavdaroglu et al. (2021) explored how a quality-based payment approach (with market prices) to incentivize farmers to improve quality. They found that this approach, combined with crop insurance, could achieve meaningful gains. Lejarza and Baldea (2022) proposed an integrated optimal decision-making framework to address the dynamic quality evolution of perishable products. Berger et al. (2022) developed an epidemiology-based framework to investigate the effects of network structure on quality controls. They proposed that aligning the network structure by focusing on the upstream-centric goods flow can decrease vulnerability to quality problems. Sajeesh et al. (2022) used a duopolistic competition model to explore how customers' awareness of health and taxation change influence firms' quality decisions. They showed that greater awareness of health and taxation could enlarge and reduce the quality gap between alternatives, respectively. However, none of the above works considered adulteration in food production. In this paper, we analyze the producers' endogenous adulteration behavior combined with quality uncertainty in raw agricultural product production, which substantially complicates the analysis.

Another stream of literature closely related to ours is on socially responsible operations in quality management, which examines opportunistic or unethical supplier behavior. Babich and Tang (2012) analyzed three mechanisms, i.e., deferred payment, inspection, and a combined mechanism, to combat adulteration. They showed that the deferred payment (inspection) mechanism can (cannot) completely eliminate the supplier's adulteration behavior, and the combined mechanism is redundant. Mu et al. (2014, 2016) studied milk farmers' deliberate adulteration behavior due to high testing costs, harmful competition, and free-riding among farmers. They proposed two interventions to make all the farmers provide high-quality milk and each competing station conduct only one mixed test and no further testing. Dong et al. (2016) found that an inspection-based approach is more effective than an external failure-based approach in deterring adulteration in a multitiered supply chain. Cho et al., (2015) explored quality and price levers to regulate adulterates, and they showed that these strategies are ineffective when the counterfeiter is deceptive. Tibola et al. (2018) proposed that frequent sampling needs to be implemented together with a cost-effective food adulteration detection system. More recently, Brooks et al. (2021) explored the impacts of the COVID-19 pandemic and Brexit on food adulterations and their impactions to consumers and manufacturers. Different from these works that modeled the detection of low-quality or defective products as exogenous, we endogenize the probability of detecting any adulteration depending on the sellers' adulteration behavior, i.e., imperfect testing. This endogeneity is more in line with practical quality inspections for agricultural products.

The work most related to this research is that of Levi et al. (2020), who investigated farms' *preemptive* and *reactive* adulteration behavior under imperfect testing in farming supply chains. They characterized the farms' equilibrium adulteration decisions under the two scenarios and analyzed the impacts of quality uncertainty, supply chain dispersion (i.e., the number of farms), traceability, and testing sensitivity on these decisions. Our work differs from that of Levi et al. (2020) in two aspects. First, Levi et al. (2020) considered homogeneous farms that adopt symmetric adulteration decisions in farm supply chains. In this work, we assume the sellers to differ from each other in their quality. Second, Levi et al. (2020) analyzed farms' adulteration behavior from the perspective of a nonprofit-seeking manufacturer and, thus, did not examine its impact on consumer choice. In this work, we assume the sellers' adulteration behavior to affect the demand in two ways, i.e., the total demand on the EP (due to its negative externality) and the consumers' choice probability among different sellers (via QPR utility).

Recently, a growing body of literature has explored the value of advanced information technologies to ensure transparency in supply chains, e.g., blockchain and IoT-based frameworks (Majdalawieh et al., 2021), and distributed trustless and secure architectures (Casino et al., 2021). Taking blockchain as an example, Cui et al. (2019) explored how blockchain-enabled traceability affects the quality contracting equilibria under serial supply chains and parallel supply chains. They showed that traceability can reduce the risk of quality uncertainty in a serial chain, but will induce conflict and can only benefit the buyer in a parallel chain. Chod et al. (2020) developed an open-source blockchain protocol to provide supply chain transparency at sale cost-effectively. Dong et al. (2022) developed a three-tier food supply chain with multiple suppliers to explore how blockchain incentivizes different stakeholders. They found that blockchain-based traceability could benefit each member in the chain but each tier is left vulnerable to its immediate downstream buyer's strategic pricing effect. Wang et al. (2021) designed a blockchain-enabled data-sharing marketplace for a stylized supply chain to achieve secured data sharing along the chain. The above works considered reliable suppliers, while this study investigates the value of traceability in combating adulteration (in Section 5.4). In this respect, this work constitutes a novel departure from their research.

## 3. Model formulation

Consider that $n$ sellers sell a certain kind of raw agricultural product on an EP. Following Liu et al. (2019) and Liu (2022), denote $\alpha_i > 0$ the overall quality level of seller $i$'s products, and we assume $\alpha_i$ to be seller $i$'s private information based on his past production experience. In practice, such a quality evaluation method for agricultural products is recommended and widely adopted by many governmental and intergovernmental organizations, such as UNECE, OECD, and WTO. For example, UNECE has developed over 100 quality standards for a wide range of perishable products (including fresh fruit and vegetables, egg products, and meat) to promote uniform quality-control procedures[14]. Without loss of generality, we assume $\alpha_1 > \cdots > \alpha_n$. It is worth noting that in agricultural production, the realized quality of a seller's output is uncertain ex ante, which is affected by many factors, such as weather conditions and infestation of pests and diseases, throughout the growing season. We represent the uncertainty in quality by random variable $\delta$ and suppose that $\delta > 0$ and $\mathrm{E}(\delta) = 1$. Therefore, if seller $i$ does not adulterate his products, he will expect to obtain their overall quality level of $\alpha_i$. Here, we implicitly assume that the uncertainty factor is perfectly correlated for all sellers, i.e., the uncertainty factor is at the "market level" (Chintapalli and Tang, 2022). This assumption is commonly adopted in the agricultural economics and operations literature (e.g., Kazaz and Webster, 2011; Alizamir et al., 2019) to ensure tractability and is reasonable when the farmers are located in countries that have similar weather conditions and, hence, are exposed to the same sources of uncertainty.

Let $x_i \in [0, 1]$ be the amount of adulterants that seller $i$ adds to his unit product, where $x_i = 0$ and $x_i = 1$ indicate that he does not adulterate and adulterates with the maximum dosage, respectively. The decisions of $x_i = 0$ and $x_i > 0$ can be interpreted as using the adulterant within and beyond the legal limit, respectively; in the latter case, seller $i$'s products will cause harm to human health (Levi et al., 2020). This interpretation is consistent with China's current food safety law, in which adulterants are defined based on whether the residue amount of highly toxic adulterants exceeds the maximum allowable limit[15]. Moreover, driven by economic interests, sellers add harmful ingredients to fresh products in an attempt to pass food safety inspections, seemingly improving food quality. For example, melamine increases the nitrogen content, which creates the illusion of increased protein concentration and improved milk quality. In addition, in the pre-emptive EMA scenario, in order to solve the problem of crop diseases and insect pests, producers will overuse harmful pesticides to "improve" product quality, and the harmful substances remaining in the products seriously threaten human health. Therefore, we denote $h(x_i)$ the resulting "quality improvement" in seller $i$'s products, and assume $h(x_i)$ to be concavely increasing with $h(0) = 1$ and $h(1) = \bar{h} > 1$, implying that the effectiveness of adulteration to improve quality is marginally decreasing (Levi et al., 2020). For simplicity, denote $\hat{\alpha}_i = \alpha_i h(x_i)$ the nominal quality level of seller $i$'s products after adulterating.

To combat adulteration, the inspectors will conduct quality testing of product samples from $t$ ($t \in [0, n]$) randomly chosen sellers. Usually, in China, the specific proportion of sampling inspections, i.e., $t/n$, is determined by administrative departments and will be announced in advance to related stakeholders[16]. Since the testing is imperfect, the probability of detecting any seller's adulteration behavior depends on the amount of adulterants he adds. Following Levi et al. (2020), this dependency can be modeled to be linearly increasing. Taken together, the probability that inspectors detect seller $i$'s adulteration behavior is $\frac{t}{n} \cdot x_i$, which is not only related to the amount of adulterants that seller $i$ adds to his unit product, but also related to the other sellers' adulteration behavior. Generally speaking, the greater the total amount of sellers' adulterants, the greater the probability that adulterated sellers will be detected. Once caught adulterating, he will be fined $\theta p_i d_i$ by the administrative departments, where $\theta > 1$ measures the intensity of penalty, and $p_i$ and $d_i$ denote his selling price and sales volume on the EP, respectively. Modeling the monetary fines for adulteration to be proportional to sales revenues (i.e., $p_i d_i$) is motivated by China's recent food safety law[17], and such a penalty structure also captures the cost of losing future business, which is especially relevant for smallholder farmers whose basic incomes rely heavily on the sales revenues of their products. It is important to point out that, based on China's current food safety law, if the EP's any seller is caught adulterating, it will not be fined by the administrative departments as long as it strictly observes the laws and regulations concerning online food safety, including (i) providing the sellers' authentic name and contact information, (ii) reviewing their business licenses, and (iii) notifying and stopping providing services when necessary[18]. Therefore, even if some of the inspectors are entrusted by the EP (usually the third-party professional testing labs), we assume the EP to be reliable in our base models, i.e., it has no incentive to cover up the sellers' adulteration behavior. Nevertheless, we extend our base model in Section 5.2 to consider an unreliable EP and propose a penalty scheme to avoid such malpractice.

As aforementioned, the sellers' adulteration behavior will create negative externality on the EP, as such behavior will be periodically released to the public. Naturally, the more quality-conscious the consumers are or the more severe the sellers' adulteration behavior revealed, the greater the negative externality is (Chen et al., 2020). For example, the greater the amount of adulterants, the easier it is for adulterated sellers to be detected, resulting in a "bad word-of-mouth effect", with some consumers boycotting the seller concerned (e.g., Ha et al., 2022), which will reduce the overall demand of the EP as well as that of other sellers. Therefore, we denote the total expected demand on the EP as $D(x_i, x_{-i}) = 1 - \rho(x_i + x_{-i})$, where $-i$ denotes all sellers other than $i$; $x_{-i} = (x_1, \ldots, x_{i-1}, x_{i+1}, \ldots, x_m)$ represents the decisions of other sellers, and $\rho > 0$ measures the level of consumers' quality consciousness. Specifically, we assume $\rho < \frac{1}{n}$ so that the demand is nonnegative. Notably, this demand function can also be interpreted as the case where the adulteration behavior is detected by consumers, the probability of which is, of course, positively impacted by the relative amount of adulterants added and their quality consciousness. To this end, if the sellers are reliable, the total demand on the EP is normalized to be one as many operations and marketing literature did (e.g., Nagarajan and Rajagopalan, 2008; Pan and Honhon, 2012; Liu et al., 2019).

Providing the sellers' prices and nominal quality (descriptions) of their products on the EP, i.e., $(p_i, \hat{\alpha}_i \delta), i = 1, \cdots, n$, a consumer will compare the QPRs of different products and buy the ones that she believes to bring maximum utility. We assume that each consumer buys from only one seller. Next, we introduce the formal definition of consumers' QPR utility.

**DEFINITION**. (*Quality-price-ratio* utility, Xie et al., 2021) The QPR of seller $i$'s products can be defined by $V_i = \frac{\hat{\alpha}_i \delta}{p_i^\gamma}$, where $\gamma > 1$ measures the price sensitivity of consumers on the EP. Furthermore, the QPR utility of seller $i$'s products can be defined by $U_i = V_i \xi_i$, where the random variable $\xi_i$ reflects the consumers' heterogeneous tastes for seller $i$'s products.

In the choice model, an outside option is allowed, i.e., consumers have no interest in any products and leave the EP with an empty hand. Denote $U_0 = V_0 \xi_0$ the consumers' reservation utility. Similar to Xie et al. (2021) and Liu (2022), we assume that $(\xi_0, \cdots, \xi_n)$ are *i.i.d.* and follow the Frechet distribution (Type II extreme value distribution), which is equivalent to assuming a Gumbel distribution (Type I extreme value distribution) random variable that yields the well-known MNL model.

Denote $\vartheta$ ($0 < \vartheta \leq \bar{\vartheta}$) the *take rate* (also referred to as *commission rate*) the EP charges from its sellers, i.e., the percentage of a seller's sales revenues paid to the EP, where $\bar{\vartheta} < 1$. For example, in the food sector, China's leading B2C online retail platforms *Tmall.com* and *JD.com* set their respective take rates between 1 to 2 percent and 3 to 5 percent[19]. Finally, we can express the sellers' and the EP's respective expected profits as follows:

$$\pi_i(x_i, p_i | x_{-i}, p_{-i}) = \mathrm{E}_\delta \left\{ p_i \cdot d_i(x_i, p_i | x_{-i}, p_{-i}) \cdot \left(1 - \vartheta - \theta \cdot \frac{t}{n} \cdot x_i\right) \right\}. \quad (1)$$

$$\pi_{ep}\left(\vartheta | \theta, \frac{t}{n}\right) = \vartheta \cdot \sum_{i=1}^{n} \mathrm{E}_\delta \{p_i \cdot d_i(x_i, p_i | x_{-i}, p_{-i})\}. \quad (2)$$

The timeline of events is as follows (see Figure 1). For *preemptive* EMA, (i) after learning about the adulteration penalty and the proportion of sampling inspections to the sellers, i.e., $(\theta, t/n)$, the EP decides the take rate $\vartheta$. (ii) Each seller simultaneously and individually decides the amount of adulterants to add to his unit product, i.e., $x_i$, by expecting the overall nominal quality level of his products, i.e., $\hat{\alpha}_i$, he will obtain. (iii) Quality uncertainty is realized, i.e., $\hat{\alpha}_i \varrho$, where $\varrho$ is any realization of $\delta$, and each seller simultaneously and individually decides his selling price, i.e., $p_i$. (iv) In the selling season, each consumer chooses to buy the products that maximize her PPR utility, and the inspectors conduct sampling inspections to the sellers' products and release the inspection results to the public. Particularly, we assume that the sales and sampling inspections activities have been in place for a sufficiently long period of time so that possible impacts of their initial history on demand are phased out. The objectives of the EP and the sellers are to maximize their own expected profits. Table 1 summarizes the key notations used throughout this paper. In the following, we denote the equilibrium values by an asterisk "∗".

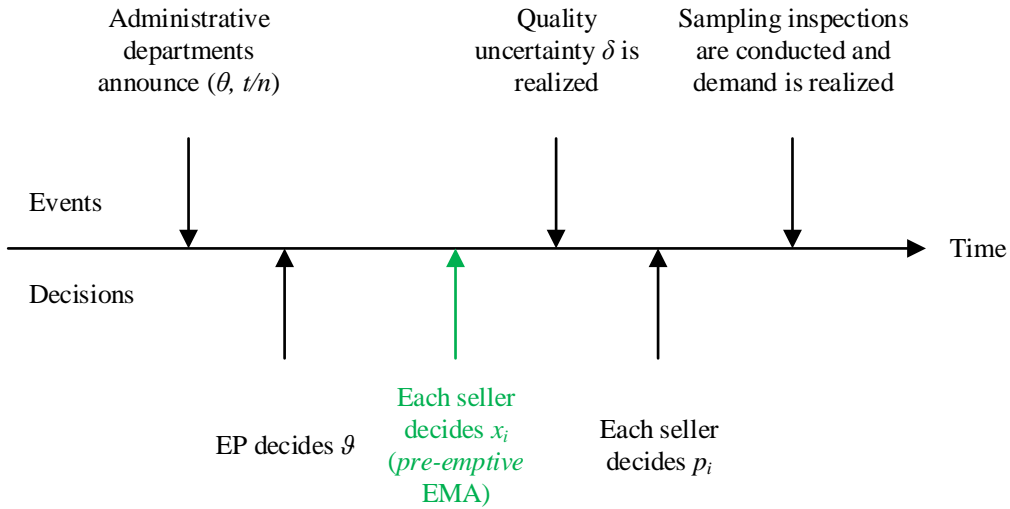

Figure 1. Timeline of events

Table 1. Key notations

| Variable | Definition |
|---|---|
| $n$ | Total number of sellers on the EP |
| $t$ | Number of farms that are randomly chosen for sampling inspections, $t \in [0, n]$ |
| $\theta$ | The intensity of penalty to sales revenues if caught adulterating, $\theta > 1$ |
| $x_i$ | The amount of adulterants seller $i$ adds to his unit product, $x_i \in [0,1]$ |
| $h(\cdot)$ | Quality improvement after adulterating, $h(\cdot) \geq 1$ |
| $\alpha_i$ | The original quality level of seller $i$'s products, $\alpha_i > 0$ |
| $\hat{\alpha}_i$ | The nominal quality level of seller $i$'s products after adulterating, $\hat{\alpha}_i = \alpha_i h(x_i)$ |
| $\delta$ | Random variable, which measures quality uncertainty, $\delta > 0$ and $E(\delta) = 1$ |
| $\rho$ | The level of consumers' quality consciousness, $0 < \rho < 1/n$ |
| $p_i$ | Seller $i$'s selling price |
| $\gamma$ | Price sensitivity of consumers on the EP, $\gamma > 1$ |
| $\vartheta$ | EP's take rate, $0 < \vartheta \leq \bar{\vartheta}$ |

## 4. Analysis

As shown in Figure 1, the EP and the sellers engage in a three-stage Stackelberg game in which the EP is the leader and the sellers are the followers. Before investigating their equilibrium decisions, we first analyze the consumers' choice probabilities among different sellers. Let $S_i$ be each individual consumer's probability of choosing to buy from seller $i$. Lemma 1 states that the choice probability for any seller is decided by and is increasing in the QPR of his products. Specifically, if the QPR of a seller's products is very high compared to the others', he can occupy a relatively large market share on the EP. It is worth noting that the sum of the consumers' choice probabilities for different sellers is smaller than 1, i.e., $\sum_{i=1}^{n} S_i < 1$, allowing for the consumers' outside option, i.e., $V_0$.

**Lemma 1**. *If the $(\xi_0, \cdots, \xi_n)$ in QPR utility functions are i.i.d. and follow the Frechet distribution, for any quality uncertainty realization $\varrho$ of the sellers' products, the consumers' choice probability of buying from producer $i$ is $S_i = \frac{\varrho V_i}{V_0 + \varrho \cdot \sum_{j=1}^{n} V_j}$.*

As Xie et al. (2021) stated, QPR-utility-based choice models have been widely adopted in the literature as the *multiplicative competitive interaction* demand model (e.g., Bernstein and Federgruen, 2004; Gallego and Wang, 2014), which is a special case of the *attraction choice model* in Gallego et al. (2006). By using the result in Lemma 1 and by (1), for any quality uncertainty realization $\varrho$, we can represent seller $i$'s expected profit on the EP as

$$\pi_i(x_i, p_i | x_{-i}, p_{-i}, \varrho) = p_i \cdot \frac{\varrho V_i}{V_0 + \varrho \cdot \sum_{j=1}^{n} V_j} \cdot [1 - \rho(x_i + x_{-i})] \cdot \left(1 - \vartheta - \theta \cdot \frac{t}{n} \cdot x_i\right). \quad (3)$$

To proceed with the analysis, we make the following assumptions:

**Assumption 1**. (i) $n \cdot \frac{\gamma - 1}{\gamma} < 1$; (ii) $\bar{h}^{\frac{1}{\gamma}} > \frac{1}{1 - n\rho}$.

Assumption 1(i) resembles Assumption 1 in Xie et al. (2021), which indicates that consumers are reasonably price-sensitive. In particular, Xie et al. (2021) have shown that if the consumers are highly price-sensitive, i.e., $n \cdot \frac{\gamma - 1}{\gamma} \geq 1$, the players will achieve a *price war* equilibrium. That is, each seller will set his selling price at the marginal cost, and no profit can be obtained. Assumption 1(ii) implies that if a seller decides to adulterate with his maximum dosage, he can substantially improve the nominal quality of his products.

Theorem 1 below characterizes the sellers' pricing game equilibrium.

**Theorem 1 (Pricing game)**. *There exists a unique equilibrium for the sellers' simultaneous pricing game. More specifically, for any quality uncertainty realization $\varrho$, seller $i$'s equilibrium selling price is $p_i^* = \left(\frac{\hat{\alpha}_i \varrho}{V^*}\right)^{1/\gamma}$, $i = 1, \cdots, n$, where $V^* = \frac{\gamma - 1}{n - (n-1)\gamma} V_0$. Consequently, the sellers will evenly share the market on the EP, that is, $S_i^* = S^* = \frac{\gamma - 1}{\gamma}$, $\forall i = 1, \cdots, n$.*

Theorem 1 shows that a seller's equilibrium pricing decision is to make the QPR of his products as competitive as the others. The intuition is as follows. Since the consumers' seller choices are determined by the QPRs of different sellers' products (see Lemma 1), a (fake) high-quality seller can take advantage of quality superiority to set a high selling price for his products without significantly sacrificing the demand. On the other hand, a low-quality seller's products are less attractive to quality-conscious consumers, so he will decrease the selling price (which is equivalent to increasing the QPR) to ensure plenty of demand on the EP, i.e., the logic of small profits but quick turnover. Consequently, the sellers' respective optimal selling prices are achieved at the QPR equilibrium, i.e., $V^*$. We can easily deduce from Theorem 1 that if the sellers are reliable, they will also evenly share the market on the EP, which means their adulteration behavior will not cause consumers' adverse selection (Stiglitz and Weiss, 1981; Liu, 2022) on the EP. However, such behavior will force consumers to afford higher prices for unqualified products that may pose threats to their health. In particular, the less unethical the sellers are (i.e., the more they adulterate), the more profits they can enjoy. Therefore, the sellers' misconduct will seriously disrupt market order on the EP.

The observations from Theorem 1 are also partially supported by data from *Taobao.com*. More specifically, we collected the monthly selling prices and sales volume data of three types of fruits, i.e., Luochuan apple, Red Fuji apple, and Zhao County pear, from *www.taobao.com* in December 2021 and calculated the coefficient of variations (CVs) of different sellers' prices and market shares for each type of fruit. Table 2 summarizes the results. We can see from Table 2 that, for all three types of fruits, their prices' CVs are obviously greater than the corresponding market shares', indicating that the sellers' very different prices lead to their relatively balanced market shares on the EP.

Table 2. CVs of price and sales volume data of three types of fruits on *Taobao.com*[20]

| Fruit | Luochuan apple | Fuji apple | Zhao County pear |
|---|---|---|---|
| Number of sellers[†] | 31 | 24 | 10 |
| CV of prices | 5.104 | 1.739 | 1.455 |
| CV of market shares | 0.187 | 0.052 | 0.099 |

[†] We only include the sellers whose monthly sales volumes are larger than 10.

Next, we analyze the sellers' equilibrium adulteration decisions. When making these decisions, quality uncertainty is yet to be realized under the *preemptive* EMA scenario. By Theorem 1 and by using $E(\delta) = 1$, we can express producer $i$'s respective expected profits under this scenario as

$$\pi_i^p(x_i, x_{-i}|p_i^*(x_i, x_{-i})) = E_\delta \left\{ \frac{\gamma-1}{\gamma} \cdot \left(\frac{\delta \hat{a}_i}{V^*}\right)^{1/\gamma} \cdot [1 - \rho(x_i + x_{-i})] \cdot \left(1 - \vartheta - \theta \cdot \frac{t}{n} \cdot x_i\right)\right\}$$

$$= \frac{\gamma-1}{\gamma} \cdot \left(\frac{\hat{a}_i}{V^*}\right)^{1/\gamma} \cdot [1 - \rho(x_i + x_{-i})] \cdot \left(1 - \vartheta - \theta \cdot \frac{t}{n} \cdot x_i\right), \quad (4)$$

Specially, we use the superscripts "$p$" and "$r$" to demote this scenario.

Before introducing an important structural property of the sellers' expected profit functions with respect to their adulteration decisions, we make the following assumption for technical convenience.

**Assumption 2.** $\frac{h'(x)}{\gamma h(x)} > \frac{n\rho}{1-n\rho x}$.

We can verify that Assumption 2 is equivalent to $\frac{\partial [h(x)^{1/\gamma} \cdot (1-n\rho x)]}{\partial x} > 0$, which means that the sellers' adulteration behavior can effectively increase their expected sales revenues after considering its negative externality. This is not surprising, as we can see from Theorem 1, that the unethical sellers will quickly increase the selling prices for their fake higher-quality products. Therefore, Assumption 2 guarantees the sellers' enough motive to adulterate their products, especially when the consumers are highly price and quality sensitive, i.e., large $\gamma$ and $\rho$, and when the competition among the sellers is fierce on the EP, i.e., large $n$.

**Proposition 1.** *For any given $x_{-i} \geq 0$, $\pi_i^p(x_i, x_{-i}|p_i^*(x_i, x_{-i}))$ and $\pi_i^r(x_i, x_{-i}|p_i^*(x_i, x_{-i}), \varrho)$ are concave in $x_i$.*

Proposition 1 indicates that the marginal value of adulteration brought to the sellers is decreasing, which can be explained as follows. As a seller adulterates more, his sales revenue will increase, but at an increasingly slow rate due to the marginally decreasing effectiveness of his adulterants. In the meantime, the detection risk of his adulteration behavior will quickly increase, and this situation could be even worse if the other sellers adulterate, further depressing the demand. As a result, the potential losses of a seller's more aggressive adulteration behavior will gradually catch up and exceed the expected revenue gain it generates.

***Remark***. As defined by EMA, fraudulent activities by sellers to substitute or add a substance to a product can increase the apparent value or reduce production costs. That is, adulteration costs are lower than the costs of standard production. We only consider sellers' equilibrium pricing decisions and adulteration decisions. In addition, Proposition 1 enables us to easily extend our models to consider these costs for the case of a linear or quadratic cost function (e.g., $a \cdot cost + b$ or $\frac{1}{2} h \cdot cost^2$) because the concave property of the sellers' expected profit functions still holds. Specifically, a quadratic cost function is widely used in agricultural models (e.g., Alizamir et al., 2019; Guda et al., 2021) to capture the diseconomies of scale prevalent in developing economies. Therefore, in our models, we do not explicitly consider the adulteration costs.

Contrary to the expectation that a lower-quality seller will adulterate more, the next theorem shows that under the *preemptive* EMA scenario the sellers' equilibrium adulteration decisions are symmetric, that is, the sellers' equilibrium strategy $x^*$ has the same structure for all $i$.

**Theorem 2 (Adulteration decisions).** *Under preemptive EMA scenario, the sellers' adulteration decisions are symmetric as follows:*

*(1) if $\theta \cdot \frac{t}{n} \geq \tau_1(1-\vartheta)$, $x^* = 0$;*

*(2) if $\theta \cdot \frac{t}{n} \leq \tau_2(1-\vartheta)$, $x^* = 1$;*

*(3) otherwise, $x^* = \hat{x}$, where $0 < \hat{x} < 1$ is decided by the following equation:*

$$\left(\frac{h'(x)}{\gamma h(x)} - \frac{\rho}{1-n\rho x}\right) \cdot \left(1 - \vartheta - \theta \cdot \frac{t}{n} \cdot x\right) = \theta \cdot \frac{t}{n}, \tag{5}$$

*where $\tau_1 = \frac{h'(0)}{\gamma} - \rho$, $\tau_2 = \left[1 + \left(\frac{h'(1)}{\gamma \bar{h}} - \frac{\rho}{1-n\rho}\right)^{-1}\right]^{-1}$, and $\tau_1 > \tau_2 > 0$.*

The symmetry of the sellers' adulteration decisions can be explained as follows. On the one hand, these sellers will share the negative externality of each other's adulteration behavior on the EP; on the other hand, their pricing decisions are ex post, i.e., they can modify their selling prices to the perceived quality to ensure the competitiveness of their products on the EP. Consequently, none of the sellers is willing to adulterate at an amount smaller than the others'. Note that $\theta \cdot \frac{t}{n}$ measures the *overall penalty risk* faced by the sellers, which takes into account the proportion of sampling inspections and the adulteration penalty. Theorem 2 shows that as this risk increases, the sellers will gradually decrease the amount of adulterants added to their products. More specifically, if the *overall penalty risk* is sufficiently high (low), they will give up adulterating (adulterate to their maximum levels); and when this risk is moderate, they will adulterate to some extent by evaluating the expected revenue gain and the potential penalty from adulteration. Levi et al. (2020) derived similar results in their Theorem 2 by assuming homogeneous players, and they showed that the *zero-adulterant* scene cannot be achieved in farming supply chains. In this paper, we extend their models to consider heterogeneous sellers, and we demonstrate that this scene can be achieved with a very high *overall penalty risk* on the EP. This finding can provide a valuable insight into food safety management in online raw agricultural product sales. That is, the administrative departments should flexibly adjust the adulteration penalty to the inspectors' limited sampling inspection capability to ensure a sufficiently high *overall penalty risk* on the EP, which is essential to combat the sellers' adulteration behavior.

Observing that the sellers' equilibrium adulteration decisions are symmetric, in the following, we drop the subscript $i$ for simplicity. The next proposition reveals how the *overall penalty risk*, the EP's take rate, and the consumers' quality consciousness impact the sellers' equilibrium adulteration decision.

**Proposition 2**. *(1) $\hat{x}$ is decreasing in $\vartheta$, $\theta$, $t$ and $\rho$;*

*(2) for any $n \geq 3$, $\hat{x}$ is increasing in $n$.*

Intuitively, the sellers will adulterate more if the *overall penalty risk* decreases or the consumers become less quality conscious. Proposition 2 confirms these intuitions. Specifically, if this risk decreases due to the smaller proportion of sampling inspections (i.e., smaller $t/n$), the sellers' *free-riding* behavior is the key driver of their more aggressive adulteration behavior, as evidence of such behavior can be found in agricultural practice (Levi et al., 2020). For example, Gadzikwa et al. (2007) provided empirical evidence of *free-riding* behavior in organic crop production, and Mu et al. (2014, 2016) provided similar evidence in the milk industry. Our intuition may also suggest that the sellers will adulterate more if the EP increases its higher take rate. However, Proposition 2 exhibits the opposite pattern. This is because the EP will not share the adulteration penalty with its sellers, as it charges a higher take rate, if a seller adulterates more to compensate his sales revenues, the increment in his possible penalty from adulteration will outpace the corresponding revenue gain. Note that, by Proposition 1, the marginal value of his adulteration behavior is decreasing and the other sellers' symmetric actions will decrease this value even further). In this case, he will try to minimize the decrease in sales revenue by reducing the possible penalty, thus adulterating less.

We now analyze the EP's take rate decision by considering its impact on the sellers' future reactions. By (2), we can express EP's expected profit as

$$\pi_{ep}\left(\vartheta|\theta,\frac{t}{n}\right) = \vartheta \cdot \frac{\gamma-1}{\gamma} \cdot \sum_{i=1}^{n}\left(\frac{\alpha_i h(x^*)}{V^*}\right)^{1/\gamma} \cdot (1 - n\rho x^*). \tag{6}$$

A higher take rate enables the EP to take a larger share of the sellers' sales revenues. However, by Proposition 2, this higher rate may decrease the sellers' sales revenues as they tend to adulterate less. To this end, how the take rate impacts EP's expected profit is determined by the trade-off between taking a larger sales revenue share and increasing the sellers' sales revenues. The next proposition investigates this impact.

**Proposition 3**. *The impact of the EP's take rate on its expected profit is as follows:*

*(1) if $\theta \cdot \frac{t}{n} \leq OPR_1, \pi_{ep}\left(\vartheta|\theta,\frac{t}{n}\right)$ is increasing in $\vartheta$;*

*(2) if $OPR_1 < \theta \cdot \frac{t}{n} < OPR_2, \pi_{ep}\left(\vartheta|\theta,\frac{t}{n}\right)$ is first increasing and then decreasing in $\vartheta$;*

*(3) if $\theta \cdot \frac{t}{n} \geq OPR_2, \pi_{ep}\left(\vartheta|\theta,\frac{t}{n}\right)$ is first increasing, then decreasing, and finally increasing in $\vartheta$, where $OPR_1 = \tau_2(1 - \bar{\vartheta})$, $OPR_2 = \tau_1(1 - \bar{\vartheta})$.*

Proposition 3 states that whether the EP can benefit from increasing the take rate depends on the *overall penalty risk* faced by its sellers. More specifically, if this risk is very low, the EP can always benefit from a higher take rate since the sellers will adulterate with their maximum dosages in this case. If this risk is moderate, the EP can benefit from increasing the take rate within a lower level; however, as this rate reaches a high level, further increasing it will decrease the EP's expected profit. This is because the sellers will first adulterate with their maximum dosages and then start to decrease their adulterants. Correspondingly, their sales revenue will first remain at a high level and then will quickly decrease (at a rate faster than the take rate increases, as we can verify from (6) that $h(\cdot)^{1/\gamma} > 1$). If this risk is very high, the EP can benefit from increasing the take rate within relatively low or high levels; otherwise, it will be worse off by increasing this rate within a median level. The reason is that, in this case, the sellers' equilibrium adulteration decision will transfer from adulterating with their maximum dosages to gradually decreasing their adulterants and finally to giving up adulteration. Correspondingly, their sales revenue will first remain at a high level, then quickly decrease, and finally remain at a low level.

In Figure 2, we use a group of numerical examples to illustrate the observations in Proposition 2 and Proposition 3. Specifically, increasing the *overall penalty risk* (by increasing $t$) and the take rate can help combat sellers' adulteration behavior. Moreover, increasing the take rate does not necessarily increase the EP's expected profit, especially when the sellers face relatively high *overall penalty risk* on the EP (e.g., $t = 20$).

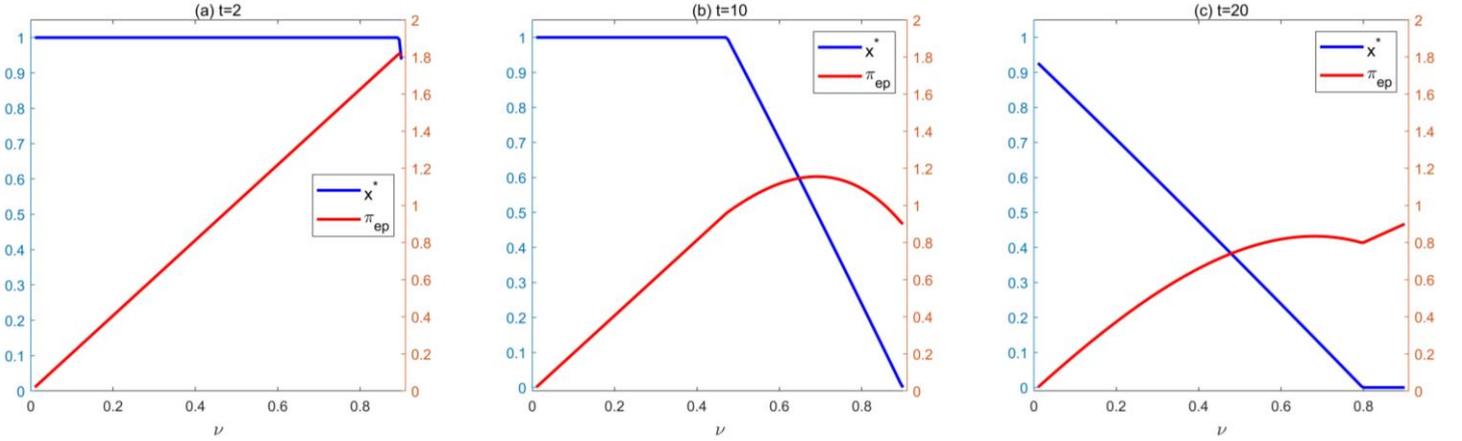

Figure 2. The impact of take rate

*Note.* Here, $n = 100, \theta = 2, \gamma = \frac{n-0.1}{n-1}, \bar{\vartheta} = 0.9, \rho = 0.0029, h(x) = 0.5x - 0.1x^2 + 1.$

Now we analyze the EP's optimal take rate decision, which is specified in Theorem 3.

**Theorem 3 (Take rate decision).** *Let $\tilde{\vartheta} = \frac{\tau_1/\tau_2 - 1}{\frac{\tau_1}{\tau_2} - \frac{1}{\bar{h}^{1/\gamma}(1-n\rho)}}$ and $OPR_3 = \tau_2\left(1 - \frac{\bar{\vartheta}}{\bar{h}^{1/\gamma}(1-n\rho)}\right)$. The EP's optimal take rate decision is given in Table 3.*

Table 3. The EP's optimal take rate decision

| $\bar{\vartheta} \leq \tilde{\vartheta}$ | | | $\bar{\vartheta} > \tilde{\vartheta}$ | | |
|---|---|---|---|---|---|
| $\theta \cdot \frac{t}{n}$ | $\vartheta^*$ | $x^*$ | $\theta \cdot \frac{t}{n}$ | $\vartheta^*$ | $x^*$ |
| $\theta \cdot \frac{t}{n} \leq OPR_1$ | $\bar{\vartheta}$ | 1 | $\theta \cdot \frac{t}{n} \leq OPR_1$ | $\bar{\vartheta}$ | 1 |
| $OPR_1 < \theta \cdot \frac{t}{n} < OPR_2$ | $1 - \frac{\theta}{\tau_2} \cdot \frac{t}{n}$ | 1 | $OPR_1 < \theta \cdot \frac{t}{n} < OPR_3$ | $1 - \frac{\theta}{\tau_2} \cdot \frac{t}{n}$ | 1 |
| $\theta \cdot \frac{t}{n} \geq OPR_2$ | $\bar{\vartheta}$ | 0 | $\theta \cdot \frac{t}{n} \geq OPR_3$ | $\bar{\vartheta}$ | 0 |

Specifically, $0 < \tilde{\vartheta} < 1$, $OPR_3 \leq OPR_2$ when $\bar{\vartheta} \leq \tilde{\vartheta}$, and $OPR_3 > OPR_2$ otherwise.

As shown in Proposition 3, the *overall penalty risk* faced by the sellers determines the take rate's impact on the EP's expected profit. Theorem 3 formalizes the link between this risk and the EP's optimal take rate decision. To be specific, if this risk is very low or high, it is optimal for the EP to set its take rate as high as possible, i.e., $\vartheta^* = \bar{\vartheta}$; if this risk is moderate, it should decrease the take rate, accordingly, to a certain level, i.e., $\vartheta^* = 1 - \frac{\theta}{\tau_2} \cdot \frac{t}{n}$. We can explain these observations as follows. In the very-low-risk case, by Proposition 3(1), the EP can always benefit from a higher take rate; and in the very-high-risk case, by Proposition 3(3), the EP can benefit from increasing the take rate to a very high level. Therefore, in both cases, the EP will set its take rate to the upper bound to maximize its expected profit. In the moderate-risk case, as the EP increases the take rate, by Proposition 3(2), its expected profit will first increase and then decrease. Therefore, its optimal take rate is achieved at the midpoint. Particularly, in this case, the higher the *overall penalty risk* the sellers face, the lower the take rate the EP will set. This finding implies that the EP will implicitly encourage sellers' adulteration behavior since increasing the higher take rate can inhibit such behavior (see Proposition 2). Again, Theorem 3 shows that the *zero-adulterant* scene can be achieved with a very high *overall penalty risk* on the EP; and in other cases, the EP's take rate decision will lead to the sellers' very severe adulteration behavior. Figure 3 illustrates these observations.

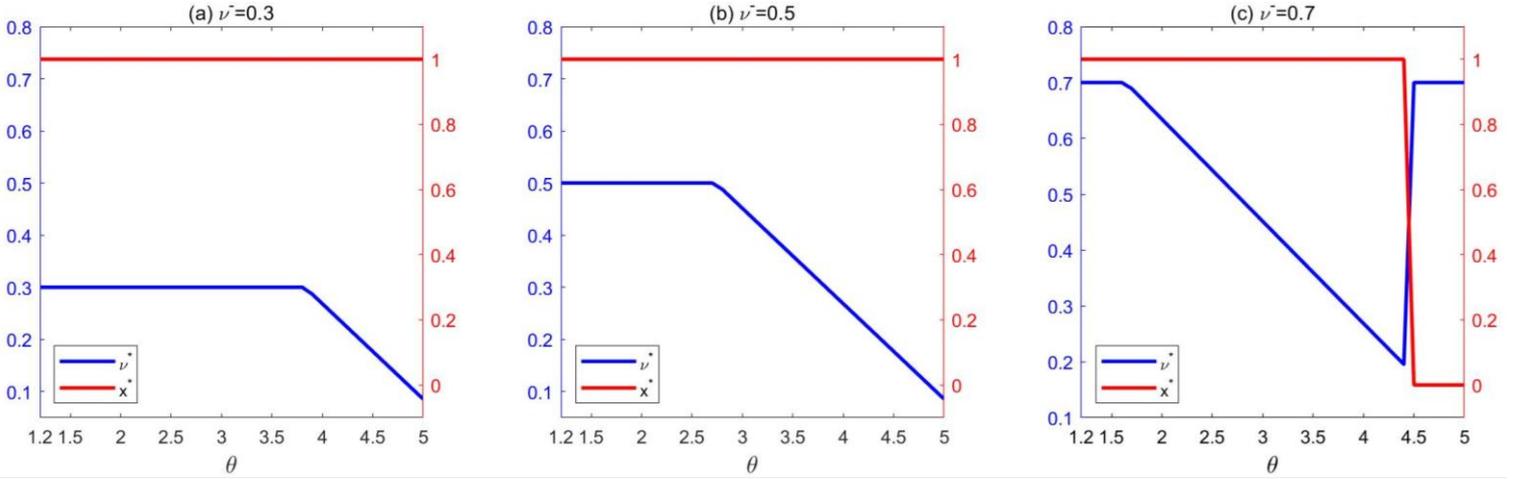

Figure 3. The optimal take rate decision

*Note.* Here, $n = 100, t = 5, \gamma = \frac{n-0.1}{n-1}, \rho = 0.0021, h(x) = 0.75x - 0.05x^2 + 1$.

The findings of Theorem 3 provide important insight into the policy implications for administrative departments. Specifically, under the imperfect testing scenario, even if the profit-maximizing EP promises to strictly enforce the laws, regulations and standards concerning food safety, the take rate decision may indulge the sellers' adulteration behavior. This unintended implication suggests that a better-designed penalty-inspection-centered approach needs to be implemented, together with the requirement for the EP's take rate. In particular, if the online sellers face a moderate *overall penalty risk* on the EP, either due to the current law's mild adulteration penalty or limited sampling inspection capability, administrative departments should require the EP to set a sufficiently large take rate from its sellers, which can help eliminate their misconduct.

## 5. Discussions

In fact, there is another scenario of EMA called reactive EMA, which occurs after the quality uncertainty is resolved and seeks to increase the perceived quality of low-quality products and create fake high-quality ones, e.g., the 2008 "melamine adulteration scandal" in China and the 2020 "honey adulteration scandal" in India. In Section 5.1, we extend our moedel to derive selles' equilibrium adulteration decisions under the *reactive* EMA scenario.

### 5.1 The reactive EMA

As shown in Figure 4, for a *reactive* EMA, the key difference is that the adulteration decisions are made after the quality uncertainty is realized.

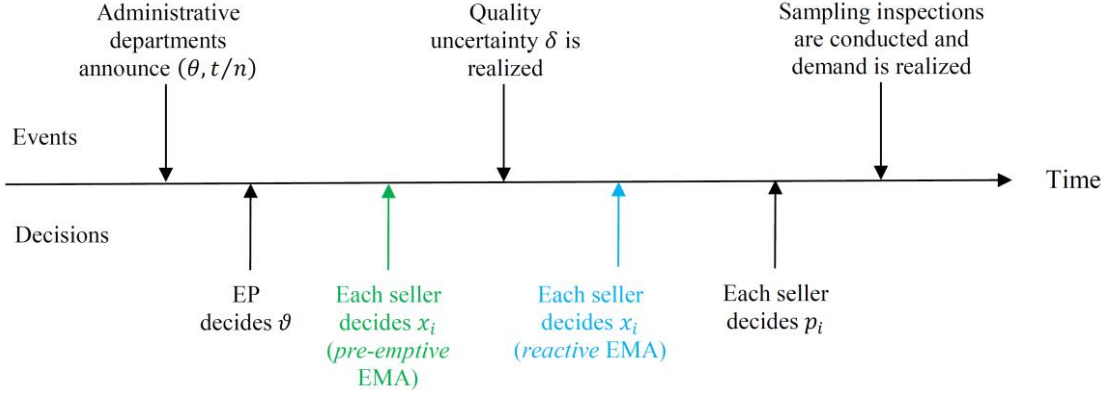

Figure 4. Timeline of the *pre-emptive* and *reactive* EMA scenarios' events

By Theorem 1 and by using $E(\delta)=1$, we can express producer $i$'s respective expected profit under this scenario as

$$\pi_i^r(x_i, x_{-i}|p_i^*(x_i, x_{-i}), \varrho) = \frac{\gamma-1}{\gamma} \cdot \left(\frac{\varrho \hat{\alpha}_i}{V^*}\right)^{1/\gamma} \cdot [1 - \rho(x_i + x_{-i})] \cdot \left(1 - \vartheta - \theta \cdot \frac{t}{n} \cdot x_i\right). \quad (7)$$

Specially, we use the superscript "r" to demote this scenario.

Similar to *preemptive* EMA, Theorem 2 also states that seller $i$'s equilibrium adulteration decision is symmetric under *reactive* EMA.

In the previous sections, we assume the EP to be reliable, i.e., they will not cover up the sellers' adulteration behavior if they are caught adulterating by the inspectors entrusted by the EP (hereinafter referred to as the EP's inspectors). However, as we have shown that the EP may indulge their adulteration behavior under certain conditions, there exists the risk that the EP turns a blind eye to their misconduct. In Section 5.2, we analyze how administrative departments can design a penalty scheme to avoid the EP's malpractice in this regard. Moreover, China's recent food safety law underlines that *when an EP promises a consumer compensation that is higher than the amount provided under the law, then the EP will have to pay the promised amount*[21]. Therefore, in Section 5.3, we investigate whether the EP has the incentive to claim a higher-than-law-requires adulteration penalty to consumers. Finally, considering the fact that many EPs have adopted advanced traceability technologies to ensure the quality transparency of their food (e.g., *JD.com*'s Blockchain Anti-counterfeiting and Traceability Platform[22]), in Section 5.4, we explore the value of traceability in eliminating adulteration on the EP.

**5.2 A penalty scheme to avoid the EP's malpractice**

Denote $t_e$ and $t_a$ ($t_e, t_a \in [0, n]$) as the number of sellers randomly chosen for quality inspections by the EP and the administrative departments' inspectors, respectively. Naturally, the EP can only manipulate the inspection results conducted by the former. Specifically, to highlight the EP's malpractice regarding food safety management, we assume that it will cover up sellers' adulteration behavior if such behavior is detected by its own inspectors. As a result, only under the case where the two types of inspectors coincidentally inspect and detect the same producer's misconduct (e.g., producer $i$, the probability of which is $\frac{t_e t_a}{n^2} \cdot x_i^2$), the administrative departments can discover the EP's malpractice by reviewing its quality inspection records. This setting is identical to the actual situation in China, as the current food safety law requires the EP to provide the necessary surveillance and monitoring information about the sellers that are caught adulterating[18]. In this case, we assume that the EP will be fined an administrative penalty of $AP$. By (6), we can derive the unreliable EP's expected profit as

$$\pi_{ep}^u\left(\vartheta|\theta,\frac{t}{n}\right) = \frac{\gamma-1}{\gamma} \cdot \sum_{i=1}^n \left(\frac{\alpha_i}{V^*}\right)^{1/\gamma} \cdot h(x^*(t_a))^{1/\gamma} \cdot \vartheta \cdot \left(1 - n\rho x^*(t_a)\right) - \frac{t_e t_a}{n} \cdot x^*(t_a)^2 \cdot AP.$$
(8)

**Theorem 4.** *To keep EP from covering up the sellers' adulteration behavior, the administrative departments should set the AP as follows:*

Table 4. The optimal setting of *AP*

| Scenarios | | $AP$ |
|---|---|---|
| 1 | $\theta \cdot \frac{t_e+t_a}{n} \leq \widetilde{OPR_1}$ | Not necessary |
| 2 | $\theta \cdot \frac{t_a}{n} \leq \widetilde{OPR_1} < \theta \cdot \frac{t_e+t_a}{n} < \widetilde{OPR_2}$ | $AP \geq AP_1$ |
| 3 | $\theta \cdot \frac{t_a}{n} \leq \widetilde{OPR_1} < \widetilde{OPR_2} \leq \theta \cdot \frac{t_e+t_a}{n}$ | $AP \geq AP_2$ |
| 4 | $\widetilde{OPR_1} < \theta \cdot \frac{t_a}{n} < \theta \cdot \frac{t_e+t_a}{n} < \widetilde{OPR_2}$ | $AP \geq AP_3$ |
| 5 | $\widetilde{OPR_1} < \theta \cdot \frac{t_a}{n} < \widetilde{OPR_2} \leq \theta \cdot \frac{t_e+t_a}{n}$ | $AP \geq AP_4$ |
| 6 | $\theta \cdot \frac{t_a}{n} \geq \widetilde{OPR_2}$ | Not necessary |

where $\widetilde{OPR_1} = \tau_2(1-\vartheta)$, $\widetilde{OPR_2} = \tau_1(1-\vartheta)$, $AP_i > 0, i = 1,2,3,4$, $AP_2 > AP_1$ and $AP_4 > AP_3$.‡

We can make the following observations from Theorem 4:

(1) If the proportion of sampling inspections from the administrative departments' inspectors is very small (Scenario 1), an $AP$ is unnecessary. This is because the sellers will severely adulterate in this case, and the administrative departments' inspectors can easily detect their misconduct and, thus, the EP's malpractice regarding food safety management. Therefore, the EP has no motive to cover up its sellers' adulteration behavior.

(2) If the proportion of sampling inspections from the administrative departments' inspectors is relatively small (Scenarios 2 and 3) or moderate (Scenarios 4 and 5), the EP has the motive to cover up the sellers' adulteration behavior. Therefore, an $AP$ is necessary. Particularly, we can see that $AP_2 > AP_1$ and $AP_4 > AP_3$, which means the EP's such motive is even stronger when its own inspectors inspect more sellers, i.e., $t_e$ increases. This is because it is more difficult for administrative departments to detect sellers' adulteration behavior and, thus, the EP's malpractice, as they will adulterate less in these scenarios than in Scenario 1.

(3) If the proportion of sampling inspections from the administrative departments'

---
‡ The specific expressions of $AP_i > 0, i = 1,2,3,4$, can be found in the proof of Theorem 4 in Appendix B.

inspectors is very large (Scenario 6), an *AP* is unnecessary as well, since the sellers will give up adulterating due to the very high detection risk.

Theorem 4 warrants that relying solely on the EP's inspectors to deter sellers' adulteration behavior is risky. We advocate the administrative departments to increase independent quality inspections to online sellers and, in the meantime, set an administrative penalty to constrain the EP's food safety management practices.

### 5.3 Whether to claim a higher-than-law-requires penalty

Increasing the adulteration penalty can help deter the sellers' adulteration behavior, and the sellers' resulting smaller sales revenues do not necessarily decrease the EP's profit, as it will optimally adjust its take rate to this penalty (see Table 3). Specifically, it will set a relatively high take rate under the case of a small or large adulteration penalty. As a result, whether the EP has the incentive to claim a higher-than-law-requires adulteration penalty again follows from a trade-off between taking a larger sales revenue share and increasing the sellers' sales revenues. Theorem 5 characterizes the condition under which the EP has the incentive to do so.

**Theorem 5.** *Only under the case where $\bar{\vartheta} \leq \tilde{\vartheta}$ and $\tilde{\theta} < \theta < \frac{OPR_2}{t/n}$, the EP has the incentive to increase $\theta$ to $\frac{OPR_2}{t/n}$; otherwise, it has no incentive to do so, where $\tilde{\theta} = \frac{\tau_2}{t/n} \cdot \left(1 - \frac{\bar{\vartheta}}{h^{1/\gamma} \cdot (1-n\rho)}\right)$.*

Theorem 5 implies that if the EP is in a relatively weak position in the market (which is reflected by its limited take rate, i.e., $\bar{\vartheta} \leq \tilde{\vartheta}$) and the law sets a moderate adulteration penalty, it has the incentive to increase this penalty. The reason is as follows. If the law sets a very low (high) penalty, the sellers will aggressively (will not) adulterate their products, and the EP will increase its take rate to the upper bound to maximize profit (see Table 3). Therefore, claiming a higher-than-law-requires penalty does no good for (is unnecessary for) the EP in the former (latter) case. If the EP dominates the market and is able to set a very high take rate, i.e., $\bar{\vartheta} > \tilde{\vartheta}$, it can ensure sufficient profit by flexibly modifying this rate to the adulteration penalty. Therefore, the EP has no incentive to do so in this case either. Otherwise, i.e., the EP's take rate is limited by its weak market position and the law sets a moderate adulteration penalty, the sellers' mild adulteration behavior will not significantly increase their sales revenues, and the EP will set a relatively small take rate in this case, i.e., $\vartheta^* = 1 - \frac{\theta}{\tau_2} \cdot \frac{t}{n} < \bar{\vartheta}$. By increasing the penalty to satisfy $\theta \cdot \frac{t}{n} \geq OPR_2$, the EP will adjust its take rate to $\vartheta^* = \bar{\vartheta}$ accordingly and can enjoy an even higher profit. In practice, consumers will regard the EP's higher-than-law-requires adulteration penalty claim as a "good-quality signal", and such a signal can further increase the EP's welfare.

The insights from Theorem 5 call for the need to set a high enough adulteration penalty and to conduct an extensive antitrust review of the market power of EPs; otherwise, these platforms can "legally" make use of the sellers' adulteration behavior to increase their profits. Indeed, China's recent food safety law has increased the penalty for producing any substandard food product from up to 10 times its value to up to 30 times[17]. Moreover, China has been curtailing monopolistic behavior in the platform economy in the last two years. For example, in 2021, China's State Administration for Market Regulation imposed a 3.44 billion yuan fine on *Meituan* (the largest online food delivery platform in China), abusing its dominant market position through its "pick one from two" practice and setting high and discriminatory take rates from its sellers, amounting to 3% of its 2020 domestic revenue[23,24].

**5.4 The value of traceability in combating EMA**

In this subsection, we explore the value of traceability in EMA elimination, as well as its implications for different stakeholders. Usually, the traceability system is developed by the EP, and its sellers will be charged the usage fee. Figure 2 gives an example of the blockchain-enabled anti-counterfeiting and traceability platform developed by *JD.com*. In some cases, the sellers' large-scale use of such a system can create a new profit stream for the EP. Denote $C_e$ as the EP's R&D expense and $C_s$ as a seller's usage fee for the tradability system. As we have shown that the sellers' equilibrium decisions under the *preemptive* and *reactive* EMA scenarios are essentially the same, in this subsection, we only consider the former scenario for simplicity. Specifically, we assume that all sellers will use the traceability service on the EP. Furthermore, we assume that the traceability system can enable credible and incorruptible traceability of the entire agricultural supply chain, so that the sellers' adulteration behavior is eliminated on a traceable EP.

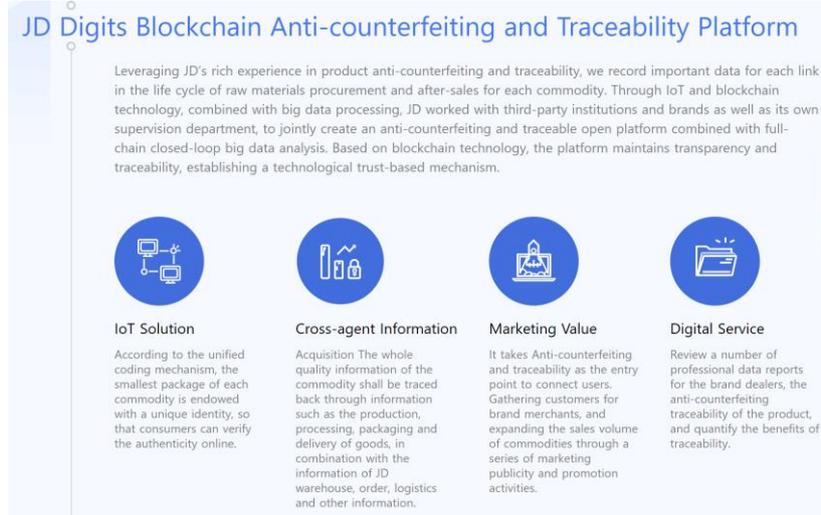

Figure 5. *JD.com*'s Blockchain Anti-counterfeiting and Traceability Platform[22]

By (2) and (3), we can derive the traceable EP's and its sellers' expected profits as

$$\pi_i^t(p_i^t|p_{-i}^t, \varrho) = p_i^t \cdot \frac{\varrho V_i^t}{V_0 + \varrho \cdot \Sigma_{j=1}^n V_j^t} \cdot (1 - \vartheta) - C_s, \tag{9}$$

$$\pi_{ep}^t(\vartheta) = \vartheta \cdot \sum_{i=1}^n p_i^t \cdot \frac{\varrho V_i^t}{V_0 + \varrho \cdot \sum_{j=1}^n V_j^t} + nC_s - C_e, \tag{10}$$

where $V_i^t = \alpha_i / p_i^{t\gamma}$. Each profit-maximizing seller must decide his selling price by considering the best response of the other sellers on the EP.

Lemma 2 derives the sellers' equilibrium pricing decisions on the traceable EP.

**Lemma 2.** *On a traceable EP, seller $i$'s equilibrium selling price is $p_i^{t*} = \left(\frac{\alpha_i}{V^*}\right)^{1/\gamma}$, $i = 1, \cdots, n$, where $V^*$ is defined the same as that in Theorem 1. Moreover, the sellers' equilibrium market shares are symmetric on the EP, i.e., $S_i^* = S^* = \frac{\gamma - 1}{\gamma}$, $\forall i = 1, \cdots, n$.*

By Lemma 2, the two stakeholders' expected profits can be rewritten as

$$\pi_i^t(p_i^{t*}) = \frac{\gamma - 1}{\gamma} \cdot \left(\frac{\alpha_i}{V^*}\right)^{1/\gamma} \cdot (1 - \vartheta) - C_s, \tag{11}$$

$$\pi_{ep}^t(\vartheta) = \vartheta \cdot \frac{\gamma - 1}{\gamma} \cdot \sum_{i=1}^n \left(\frac{\alpha_i}{V^*}\right)^{1/\gamma} + nC_s - C_e. \tag{12}$$

Furthermore, their total expected profits can be derived as

$$\pi_{total}^t = \frac{\gamma - 1}{\gamma} \cdot \sum_{i=1}^n \left(\frac{\alpha_i}{V^*}\right)^{1/\gamma} - C_e. \tag{13}$$

Theorem 6 specifies the conditions under which the traceability system can add value to these two stakeholders.

**Theorem 6.** *The traceability system can add value to the two stakeholders under the following two cases:*

*(1) when $\bar{\vartheta} \leq \tilde{\vartheta}$, $C_e \leq \hat{C}_e$ and $OPR_4 < \theta \cdot \frac{t}{n} < OPR_2$;*

*(2) when $\bar{\vartheta} > \tilde{\vartheta}$, $C_e \leq \hat{C}_e$ and $OPR_4 < \theta \cdot \frac{t}{n} < OPR_3$.*

*Otherwise, this system cannot add value to them, where $\hat{C}_e = \frac{\gamma - 1}{\gamma} \cdot \sum_{i=1}^n \left(\frac{\alpha_i}{V^*}\right)^{1/\gamma} \left[1 - \bar{h}^{1/\gamma} \cdot (1 - n\rho) \cdot (1 - \theta \cdot \frac{t}{n})\right]$, $OPR_4 = 1 - \frac{1}{\bar{h}^{1/\gamma} \cdot (1 - n\rho)}$, and $0 < OPR_4 < 1$.*

Theorem 6 states that the traceability system can add value to the two stakeholders when the EP's R&D expense of this system is small enough and the sellers' *overall penalty risk* on a nontraceable EP is moderate. Intuitively, under the very high (low) risk case, adopting this traceability system is unnecessary (will substantially decrease the sellers' sales revenues) since they will not (will severely) adulterate their products on a nontraceable EP. Therefore, the two stakeholders are not willing to adopt this system in these two cases due to the extra costs (the decrease in their profits) in the former (latter) case. Under the moderate risk case, the sellers' mild adulteration behavior can only slightly increase their sales revenues on a nontraceable EP. Instead, in this case, adopting the traceability system can increase the two stakeholders' total expected profits by avoiding the sellers' adulteration penalty, which is especially the case when there are a large number of sellers on the EP.

Theorem 6 indicates that the traceability system can create value and is a good supplement to the penalty-inspection-centered approach to combat the sellers' adulteration behavior on the EP. However, such value is very limited if the adulteration penalty is small and the inspectors' sampling inspection capability is poor. Therefore, Theorem 6 again underlines the essential need to incorporate a systemic perspective to complement a penalty-inspection-centered approach to reduce EMA risk on the EP.

Under the cases where the traceability system can add value to the two stakeholders, the EP can achieve a win–win situation with its sellers by setting a reasonable usage fee for this system. In fact, this reasonable usage fee is the prerequisite for adopting the traceability system to the EP, which is specified in the following proposition.

**Proposition 4**. *Under the cases where the traceability system can add value to the two stakeholders, both of them can be better off on the traceable EP when $\check{C}_s \leq C_s \leq \hat{C}_s$, where $\check{C}_s = \frac{\vartheta^*(\gamma-1)}{n\gamma} \cdot \sum_{i=1}^{n} \left(\frac{\alpha_i}{V^*}\right)^{1/\gamma} \left(\overline{h}^{\frac{1}{\gamma}} \cdot (1-n\rho) - 1\right) + \frac{c_e}{n}$, and $\hat{C}_s = \frac{\gamma-1}{\gamma} \cdot \left(\frac{\alpha_i}{V^*}\right)^{1/\gamma} \cdot \left[n\rho(1-\vartheta^*) + (1-n\rho) \cdot \theta \cdot \frac{t}{n}\right]$.*

## 6. Modeling calibration

In this section, we used the real-world data derived from *Taobao.com* to calibrate our model parameters and examined how well our theoretical analysis aligns with empirical evidence. Specifically, we chose the Luochuan apple, a characteristic fruit of Luochuan County, Shaanxi Province, China, as the research subject. In the following, we first give an overview of the dataset.

There were 88 Luochuan apple sellers on *Taobao.com* in December 2021, and we collected their selling prices and sales volume data. An example of such data is shown in Figure 6. To rule out trivial cases, if a seller sells more than one product specification (in terms of product size and packaging, etc.), we took the average of their different unit selling prices to approximate his selling price. The sellers' sales volume data vary widely from 0 to 46700, and only 31 sellers sold more than 10 products in that month. In particular, the total sales volume of these 31 sellers takes up more than 90 percent of that of the platform. To this end, following Xie et al. (2021), we treated the other sellers' (whose sales volumes are smaller than 10) products as the outside option in the QPR-based consumer choice model. Then, we have $n = 31$. We refer readers to Appendix A for the 31 sellers' more detailed data used in the numerical studies. Furthermore, to satisfy Assumption 1(i), we set $\gamma = \frac{n-0.1}{n-1} = 1.03$. According to Xie et al. (2021), for the QPR model, given any seller's selling price $p_i$, the quality of his products can be derived by

$$\alpha_i = \frac{\frac{\gamma-1}{\gamma}}{1 - n \cdot \frac{\gamma-1}{\gamma}} \cdot p_i^{\gamma}. \tag{14}$$

Considering that Chinese EPs usually set their take rates between 1 and 5 percent in the food sector[19], we set $\vartheta = 0.05$. Moreover, based on the sampling inspections conducted by Yulin's (Guangxi Province, China) government to online sellers in 2020[16], the proportion of which was 1/3, we set $t = 12$ in our base models. Furthermore, inspired by Levi et al. (2020), we used the quadratic function $h(x) = a(bx - x^2) + 1$ in the numerical studies, and we require $a > 0$ and $b > 2$ to satisfy the concavely increasing property. Obviously, the greater $a$ and $b$ indicate the higher effectiveness of the sellers' adulteration behavior regarding (fake) quality improvement. Specifically, we set $(a, b) = \left(\frac{1}{15}, 15\right)$ in our base models. Finally, we set $\rho = 0.8 \cdot min\left\{\frac{1}{n} \cdot \left(1 - \bar{h}^{-\frac{1}{\gamma}}\right), \frac{1}{n} \cdot \frac{a(b-2)}{\gamma\bar{h}+ab}\right\}$ to simultaneously satisfy Assumption 1(ii) and Assumption 2.

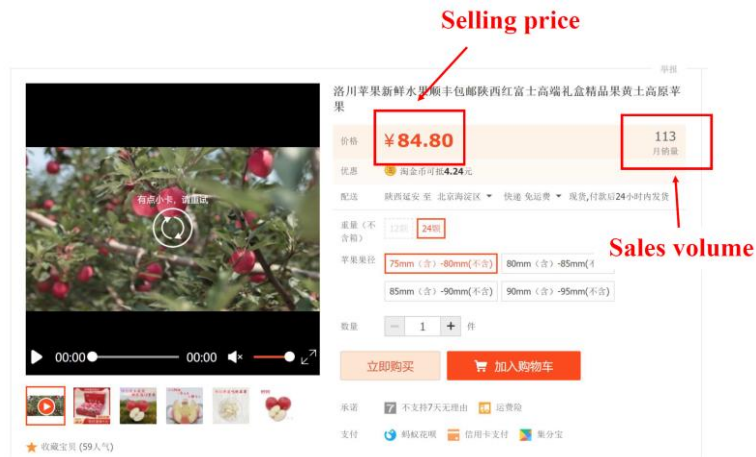

Figure 6. An example of selling prices and sales volume data

We conducted three groups of experiments with respect to $a$, i.e., $a = \frac{1}{5}, \frac{1}{10}, \frac{1}{15}$, and for each group of experiments, we varied $\theta$ from 3 to 10 with an interval of 0.5. We report the sellers' equilibrium adulteration decisions, i.e., $x^*$, and the EP's and the sellers' total expected profits, i.e., $\pi^{n*}$, under three sampling inspection scenarios, i.e., $t = 3$, $t = 8$, and $t = 12$. In addition, we report the two stakeholders' total expected profits on a traceable EP (i.e., the black dotted line $\pi^{t*}$) as the benchmark. Figure 7 depicts the results.

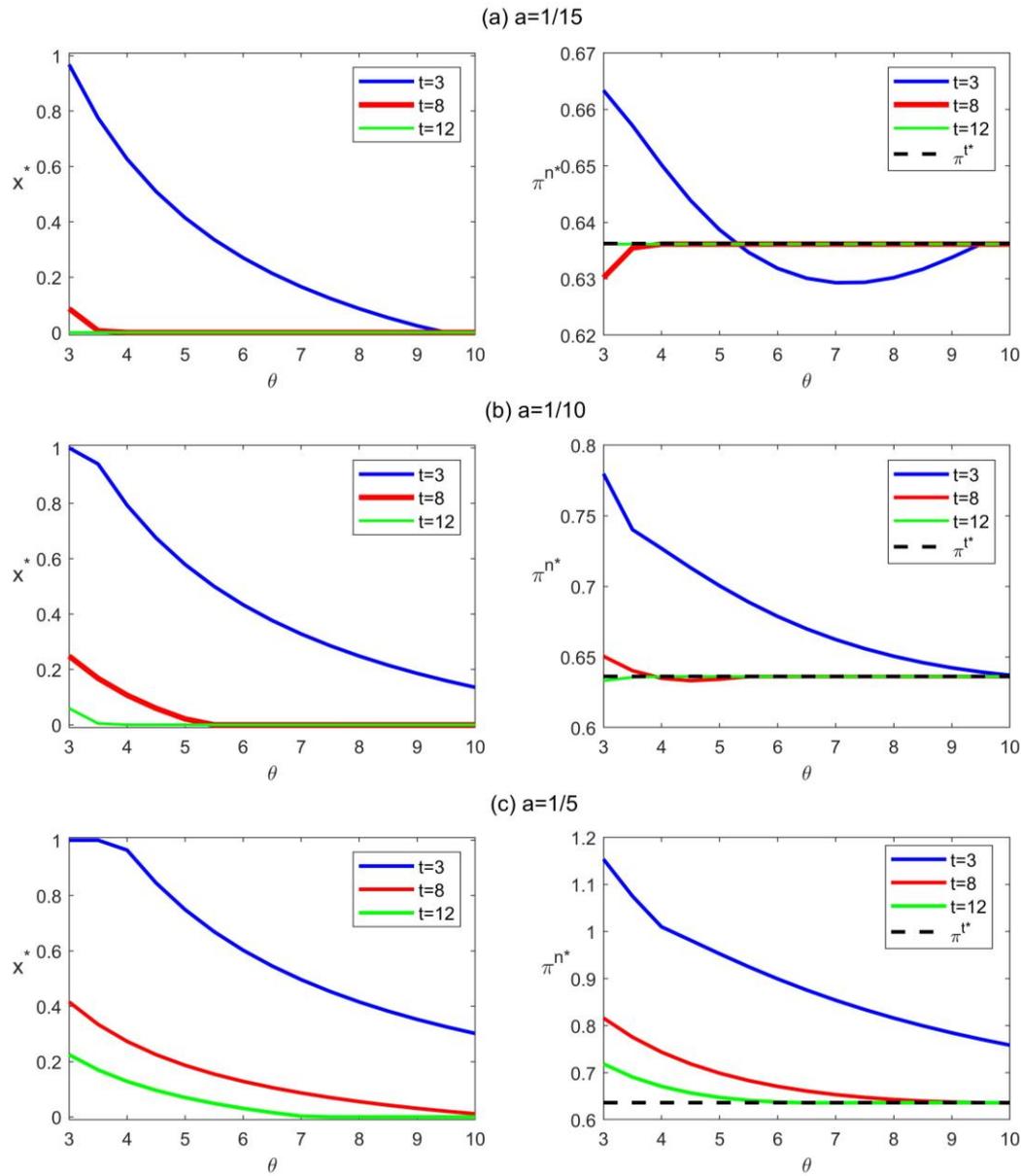

Figure 7. EMA risk analysis of *Taobao.com*'s Luochuan apple

We can make the following observations from Figure 7:

(1) Consistent with our intuition, increasing the adulteration penalty (i.e., $\theta$) or the proportion of sampling inspections (i.e., $t/n$) will deter the sellers' adulteration behavior. Moreover, if the adulterants become more effective (i.e., a larger $a$), the sellers will adulterate more aggressively. This is straightforward because they can enjoy the same high selling prices with a smaller amount of adulterants and, thus, a lower detection risk on the EP. Consequently, the greater equilibrium amount of adulterants follows from the trade-off point between the revenue gain and potential penalty moving toward the right side.

(2) Intuitively, if the *overall penalty risk* decreases, the two stakeholders' total expected profit will increase since the sellers will adulterate more (thus generate higher sales revenues) and the potential penalty decreases. To our surprise, Figure 7 shows that this intuition only holds (i) when the adulterants are very effective (e.g., $a = 1/5$) and (ii) when the adulterants are moderately effective (e.g., $a = 1/10$) while the proportion of sampling inspections is small (e.g., $t = 3$). This phenomenon can be explained as follows. The sellers will adulterate with a relatively large dosage in both cases, indicating the high detection risk of their misconduct. Therefore, if the *overall penalty risk* decreases, the two stakeholders' total expected profit will increase because of the significant decrease in the adulteration penalty.

(3) However, when the adulterants are lowly effective (e.g., $a = 1/15$), increasing this penalty will first decrease and then increase (will increase) the two stakeholders' total expected profit under the case of a small (moderate) proportion of sampling inspections, i.e., the blue (red) line in Figure 7(a). The reason is that the sellers benefit limitedly from their misconduct. In the former case, they will adulterate with a relatively large dosage, and a slight increase in the adulteration penalty will obviously increase their potential losses. Therefore, their total expected profit will decease as the sellers substantially decrease their adulterants, i.e., much smaller sales revenues. As the sellers' adulterants decrease to a very small dosage, the already very low detection risk will induce them to decrease this dosage more slowly (see the left figures in Figure 7), and their revenue gain starts to outpace the potential penalty from adulteration. Consequently, their total expected profit will start to increase. In the latter case, they will adulterate with a moderate dosage, and increasing this penalty will further decrease this dosage. As a result, the two stakeholders' total expected profit will increase due to the smaller penalty. We can see the similar (despite less obvious) patterns in Figure 7(b).

These observations offer the following insights for policy-makers and commercial entities in online raw agricultural product sales. First, we once again highlight the penalty-inspection-centered approach to deter EMA, especially when the unethical sellers can use the adulterants sophisticatedly. Second, under certain conditions, increasing the *overall penalty risk* can add value to the EP and its sellers; however, the traceability system is a better choice for the two stakeholders in these conditions, as we can observe from Figure 7(a) and (b) that the black dotted line generates the highest profit therein. This is because system adoption can help avoid the (many) sellers' adulteration penalty. Therefore, Figure 7 justifies our insights in Section 5.4 that the traceability system complements the penalty-inspection-centered approach and is a promising solution to circumvent the EMA risk on the EP.

## 7. Conclusions

In this paper, we attempt to examine online raw agricultural product sellers' *preemptive* EMA behavior under quality uncertainty and imperfect testing, and investigate how inspectors' proportion of sampling inspections, administrative departments' adulteration penalty, and EP's take rate jointly impact their behavior. Specifically, the first two factors together decide the *overall penalty risk* faced by the sellers on the EP. We develop a three-stage Stackelberg game to capture the strategic interactions among the EP, sellers, consumers, and administrative departments and proceed with the analysis in reverse order. First, based on the QPR-based consumer choice model, we analyze the sellers' selling prices. We find that, at equilibrium, a seller's pricing decision is to make the nominal QPR of his products be as competitive as others. Consequently, they will evenly share the market on the EP. Second, we characterize the sellers' equilibrium adulteration decisions. Notably, at this stage, quality uncertainty has not been realized under the *preemptive*EMA scenario. Surprisingly, we find that the heterogeneous sellers will take symmetric adulteration in both scenarios, which is determined by the *overall penalty risk* on the EP. Overall, the sellers' equilibrium amount of adulterants is decreasing in this risk, and the *zero-adulterant* (*maximum-adulterant*) scene will be achieved with a sufficiently large (small) risk. Counter to our intuition that the EP's higher take rate will induce the sellers to adulterate more, we show that this higher rate actually inhibits such behavior. To this end, third, we analyze the EP's optimal take rate decision considering its future impact on the sellers' decisions. We find that profit-maximizing EP will increase its take rate to the upper bound when the sellers face very low or high *overall penalty risk* and will decrease this rate to a certain level when they face moderate risk. In particular, the EP's decision will drive the sellers to severely adulterate under the low- and moderate-risk case, indicating that the EP may indulge the sellers' misconduct in these cases.

Next, we discuss three levers that administrative departments or the EP can use to mitigate sellers' adulteration behavior. First, we analyze how administrative departments can design an administrative penalty to deter the EP's malpractice regarding EMA elimination. We find that such a penalty is necessary (unnecessary) when the proportion of sampling inspections from the administrative departments' inspectors is relatively smaller or moderate (very small or large). Second, we investigate whether the EP has the motive to claim an adulteration penalty higher than law requires to its consumers. We show that the EP has such motive only when its take rate is limited, while the law sets a moderate penalty. Third, we exploit the value of traceability systems in deterring EMA. We find this system can add value to the EP and the sellers when its R&D expense is small enough and the sellers face a moderate *overall penalty risk* on the nontraceable EP. Generally, the above three levers can only act as supplements to the penalty-inspection-centered approach, which is essential to EMA elimination.

Finally, we carry out modeling calibration with real-world data from *Taobao.com* and investigate how the effectiveness of adulterants and *overall penalty risk* on the EP impact the sellers' adulteration behavior and the two stakeholders' total expected profit. We find that increasing the former (latter) factor will encourage (depress) sellers' adulteration behavior. Moreover, we show that the greater *overall penalty risk* does not necessarily decrease the EP's and the sellers' total expected profit. More specifically, when the sellers' adulterants are lowly effective, as the adulteration penalty increases, the two stakeholders' total expected profit will first decrease and then increase (will increase) under the case of a small (moderate) proportion of sampling inspections. However, adopting the traceability system can make both stakeholders better off than increasing this penalty in these conditions.

The problem analyzed in this paper can be further explored in a number of directions. For example, we assume the quality uncertainty to be perfectly correlated for all sellers. Naturally, sellers may exhibit different capabilities in quality control during agricultural production, and such capabilities will fundamentally impact their *preemptive* EMA behavior. It would be interesting to investigate how a seller's quality control impacts such behavior. Moreover, in Assumption 2, we assume the negative externality of the sellers' EMA behavior to be small enough so that their adulteration behavior can effectively increase their sales revenues. In practice, this negative externality could be significant due to the very easy spread of information through social media and consumers' increasing food safety concerns. Therefore, relaxing Assumption 2 would provide more useful insights to understand sellers' EMA behavior. Furthermore, we assume the total demand on the EP is stable. In reality, there are usually more than two EPs competing for demand in the market, and they often design various discount and subsidy policies to increase their orders. It would also be a fruitful direction to analyze the impact of the EP's marketing strategy on the sellers' EMA behavior. Therefore, the decisions of sellers and EPs are influenced by many factors. Although heterogeneity is taken into account, more nuanced models can be developed to understand and identify the behavior of each agent.

**Endnotes**

1. Reported by *China Daily* on December 11, 2020: E-commerce Business Liable for Food Safety.
2. This statement can be found in the report *Food Fraud and "Economically Motivated Adulteration" of Food and Food Ingredients* released by U.S. Congressional Research Service.
3. Data source, FDA: https://www.fda.gov/food/compliance-enforcement-food/economically-motivated-adulteration-food-fraud
4. The details can be found in https://www.thecasecentre.org/products/view?id=124875
5. Every March 15, the China Central Television (CCTV) hosts a two-hour prime-time show that exposes issues from f ake and defective products. The details can be found in https://www.31du.cn/news/2020-315-2.html
6. The statement can be found in the official website of Chinese government: http://www.gov.cn/gongbao/content/2007/content_764220.htm
7. Reported by *South China Morning Post* on May 5, 2020: 13 Million Farmers Are Selling Goods Online in China.
8. Data source, *Statista*: https://www.statista.com/statistics/1253521/china-online-retail-value-of-agricultural-products/
9. Data source: https://baike.baidu.com/item/%E4%BA%AC%E4%B8%9C%E7%94%9F%E9%B2%9C/22569042

10. More details on *Taobao.com*: https://rulechannel.taobao.com/?type=detail&ruleId=14&cId=1157#/rule/detail?ruleId=14&cId=1157
11. More details can be found in https://rule.jd.com/rule/ruleDetail.action?ruleId=2754
12. These statements can be found in the official website of Chinese government, in Chinese, http://www.gov.cn/gongbao/content/2017/content_5174527.htm
13. An example can be found in https://gkml.samr.gov.cn/nsjg/rzjgs/202201/t20220112_339080.html
14. The statement can be found in the official website of UNECE, https://unece.org/agricultural-quality-standards.
15. Food Safety Law of the People's Republic of China, http://www.gov.cn/zhengce/2015-04/25/content_2853643.htm.
16. An example can be found in http://www.yuhuan.gov.cn/art/2020/12/11/art_1229433488_59037502.html
17. The statement can be found in https://www.uschina.org/understanding-china%E2%80%99s-food-safety-law
18. *Measures for the Investigation and Punishment of Illegal Acts Concerning Online Food Safety* (2016): http://www.gov.cn/gongbao/content/2017/content_5174527.htm
19. Data source: https://www.tmogroup.com.cn/more/ecommerce-platform/31808/
20. The data is collected on 12/9/2021. We delete the data of those sellers that set abnormally high or low prices compared to the others, and those sellers whose sales volumes are smaller than 10 are treated as the outside option.
21. The statement can be found in https://www.china-briefing.com/news/e-commerce-platform-operators-liabilities-in-china-food-safety/.
22. More details can be found in https://blockchain.jd.com/en/.
23. Data source: https://www.scmp.com/tech/big-tech/article/3151675/china-fines-meituan-less-expected-us530-million-monopolistic
24. Data source: http://www.xinhuanet.com/tech/2021-06/16/c_1127565938.htm